\documentclass[graybox]{svmult}

\usepackage{type1cm}   
\usepackage{makeidx}        
\usepackage{graphicx}  
\usepackage{multicol}       
\usepackage[bottom]{footmisc}
\usepackage{hyperref}
\usepackage{newtxtext}      
\usepackage[varvw]{newtxmath}   

\makeindex          

\begin{document}

\title*{Generative Models for Parameter Space Reduction applied to Reduced Order Modelling}

\author{Guglielmo Padula\orcidID{0009-0005-5269-3296} and\\ Gianluigi Rozza\orcidID{0000-0002-0810-8812}}
\institute{Guglielmo Padula \at International School for Advanced Studies, Via Bonomea 265, \email{gpadula@sissa.it}
\and Gianluigi Rozza \at International School for Advanced Studies, Via Bonomea 265, \email{grozza@sissa.it}}
\maketitle
\abstract{
Solving and optimising Partial Differential Equations (PDEs) in geometrically parameterised domains often requires iterative methods, leading to high computational and time complexities.
One potential solution is to learn a direct mapping from the parameters to the PDE solution. Two prominent methods for this are Data-driven Non-Intrusive Reduced Order Models (DROMs) and Parametrised Physics Informed Neural Networks (PPINNs). However, their accuracy tends to degrade as the number of geometric parameters increases.
To address this, we propose adopting Generative Models to create new geometries, effectively reducing the number of parameters, and improving the performance of DROMs and PPINNs.
The first section briefly reviews the general theory of Generative Models and provides some examples, whereas the second focusses on their application to geometries with fixed or variable points, emphasising their integration with DROMs and PPINNs. DROMs trained on geometries generated by these models demonstrate enhanced accuracy due to reduced parameter dimensionality. For PPINNs, we introduce a methodology that leverages Generative Models to reduce the parameter dimensions and improve convergence. This approach is tested on a Poisson equation defined over deformed Stanford Bunny domains.}

\section*{Introduction}

Numerical simulations based on Partial Differential Equations (PDEs) are foundational tools across various scientific domains, allowing us to forecast behaviour of complex systems. From aeronautics \cite{shapeopt} to cardiovascular \cite{siena} applications, solving PDEs is of paramount importance, particularly when the system depends on different parameters. A special case of the latter is Geometrically Parametrised Partial Differential Equations (GPPDEs), where the physical PDE domain is parametrised. Optimising the PDE domain and then computing the associated solution is found in various scientific applications; see \cite{shapeopt, shapeopt3,ivagnespropeller,shapeopt2} for a far-from-exhaustive list of references.
 \\ As a concrete example, consider the design process of a naval ship hull: one would like to find the correct geometry that minimises the drag force of water on the hull and concurrently reduces fuel consumption.
However, solving these types of simulations is extremely costly in terms of computational time (it may take days or weeks to compute the simulations, even on High-Performance Computing infrastructures), and energy consumption (to cool down the High-Performance Computing infrastructures). \\
One solution is to employ machine learning models to directly map parameters to their associated PDE solutions. There are two ways to accomplish this, either by generating a series of snapshots with different parameters and then learning the map in a supervised manner, which is the approach of Data-Driven Non-Intrusive Reduced Order Models (DROMs) \cite{siena, gonnella,  ivagnesspatial, khamlich,  papapicco, pichi, ezyrb, sparsegrids1, sparsegrids2, gpr1, gpr2, rbf1, rbf2, pce1, pce2} or by learning the map while performing the simulation, which is the approach of Parametrised Physics-Informed Neural Networks (PPINNs) \cite{ppinn, pina}. 
\\  

DROMs and PPINNs can be combined, an example is the Digital Twin, in which it is necessary to compute the optimal policy given the current state of the physical environment, and it is crucial to obtain it in the shortest possible time. A Digital Twin can incopororate both real world or simulation data (DROMs approach), and some physics knowledge in the form of equations (PPINNs approach). In this case, it is also useful to perform Uncertainty Quantification \cite{garom,rozza}.\\

Although DROMS and PPINNS are promising in speeding up the simulation, their performance decreases as the parameter dimension increases. To overcome these problems, parameter reduction techniques have been applied \cite{shapeopt, cgm}.
However, increasing the parameter space remains a bottleneck in generating new geometries. From this perspective, Generative Models have recently been applied to create new geometries \cite{cgm,ranjan,yuan} thanks to the approximation power of Neural Networks.\\
\\
In this chapter, we provide an overview of Generative Models and study how they can be applied for parameter reduction and for designing geometries with equal or different topologies with an application to DROMs and PPINNs. \\
The chapter is structured as follows. Section \ref{sec:one} deals with the mathematics behind Generative Models. Section \ref{sec:two} introduces a technique that employs a Generative Model to generate constrained geometries, while performing parameter reduction on geometries with a fixed topology. Then, we develop an extension of the previously mentioned technique for parameter reduction to learn arbitrary topology parametrised geometries and their corresponding GPPDE solutions via PPINNs. Conclusions and remarks are presented in Section \ref{sec:3}.

\section{Introduction to Generative Models}\label{sec:one}

In this section, we provide an informal overview of Generative Models as machine learning models designed to create new data. Following this, we embark on a systematic exploration, starting with the concept of divergence. A Generative Model is a machine learning model that learns to create new data that resemble a given set of training data. This is achieved by choosing a parametric distribution based on neural networks \cite{dlbook} as they can approximate any continuous function that is zero at infinity \cite{univapprox}, and by optimising their parameters using divergences (or an approximation of them), which are measure of similarity of probability distributions. \\
Once the model is trained, it can generate new samples that resemble those in the training data.\\ These concepts will be formalised below: in the first place we will define the concept of divergences and study their properties, in the second place we will introduce the Generative Models. In the last place, we will briefly describe the architecture of the Generative Models that we will adopt later in the applications, along with the losses applied during training. \subsection{Divergences} First, we define the concept of divergence and then present three examples. With $Prob(\mathbb{R}^{n})$ we denote the space of probability distributions in $\mathbb{R}^{n}$ with the Lebesgue $\sigma$-algebra and finite first absolute moment. 
\begin{definition}[Divergence] A divergence is a function \begin{equation}D:Prob(\mathbb{R}^{n})\times Prob(\mathbb{R}^{n})\rightarrow \mathbb{R},\end{equation} such that \begin{itemize}     \item $D(P_{1},P_{2})\ge 0 \quad \forall P_{1},P_{2}\in Prob(\mathbb{R}^{n})$.      \item $D(P_{1},P_{2})= 0 \Leftrightarrow P_{1}=P_{2}.$   \end{itemize} \end{definition} An example of divergence is the Kullback-Leibler divergence.\\
It is based on the concept of the Radon-Nikodym derivative \cite{rniky}.
\begin{proposition}[Radon-Nikodym derivative]
Let $P_{1}, P_{2} \in Prob(\mathbb{R}^{n})$, such that $P_{1}\ll P_{2}$ (which means that $P_{1}$ is absolutely continuous with respect to $P_{2}$), then exists $f\in \mathbb{R}^{n}\rightarrow [0,1]$, unique up to a set of measure zero with respect to $P_{2}$, such that
\begin{equation}
P_{1}(B)=\int\limits_{B}f(x)P_{2}(dx).
\end{equation}

\end{proposition}

\begin{definition}[Kullback-Leibler (KL) divergence] Let $P_{1},P_{2}\in Prob(\mathbb{R}^{n})$ such that $P_{1} \ll P_{2}$, we define the Kullback-Leibler divergence as  \begin{equation}D_{KL}(P_{1}||P_{2})=\int\limits_{ \mathbb{R}^{n}} \log \left(\frac{P_{1}(da)}{P_{2}(da)}\right) P_{1}(da),\end{equation} where $\frac{P_{1}({d} a)}{P_{2}({d} a)} \text { is the Radon-Nikodym derivative of } P_{1} \text { with respect to } P_{2}$. \end{definition} This function is not symmetric and is employed in Variational Autoencoders (VAE) \cite{vae}, a family of Generative Models, and in classical statistics.\\

\begin{warning}{Attention}
When we declare $P$ to be absolutely continuous without specifying the base measure, the base measure is the Lebesgue measure.
\end{warning}

Indeed if $P_{1}$ and $P_{2}$ are absolutely continuous they have densities $p_{1}$, $p_{2}$, and it can be proven that
\begin{equation}
  D_{KL}(P_{1}||P_{2})=\int\limits_{\mathbb{R}^{n}} \log \left(\frac{p_{1}(a)}{p_{2}(a)}\right) p_{1}(a)da.
\end{equation}

If $p_{2}$ depends on a parameter $\theta$ (in this case we denote it as $p_{\theta}$ and the corresponding probability distribution is represented by $P_{\theta}$), and $a_{1},...,a_{m}$ are i.i.d. from $p_{1}$ because of Monte Carlo integration we have \begin{equation}
D_{KL}(P_{1}||P_{\theta})\approx -\sum \log(p_{\theta}(a_{i})),\end{equation} which is the log-likelihood function if $p_{\theta}$ represents a statistical model. In this specific case, minimising the KL divergence is equivalent to performing a maximum likelihood estimation. \\
It may happen \cite{vae,ddpm} that $P_{2}$ is defined of an "larger" domain $\mathbb{R}^{n}\times \mathbb{R}^{k}$. In this case we define as $P_2(da)$ the marginal distribution on $\mathbb{R}^{n}$, $P_2(ds)$ the marginal distribution on $\mathbb{R}^{k}$, $P_{2}(ds|A)$ the conditional distribution on $\mathbb{R}^{k}$ given $A\subset \mathbb{R}^{n}$, and $P_{2}(da|S)$ the conditional distribution on $\mathbb{R}^{n}$ given $S\subset \mathbb{R}^{k}.$
\begin{warning}{Notation for distribution and densities}
$P_{2}(ds)$, $P_{2}(da)$, $P_{2}(da,ds)$, $P_{2}(ds|A)$, $P_{2}(da|S)$ represent different distributions, and should be defined with different symbols. With a slight abuse of notation,  we will keep the same symbol to make this manuscript more readable. If we do not specify the arguments of $P_{2}$, we will intend the full distribution. 
The same applies for every other probability distribution.\\
All the distributions will be denoted with a $P$ and a subscript. 
Furthermore, when we write $p$ with a subscript, we intend a probability density associated to an absolutely continuous probability $P$ that has the same subscript. 
\end{warning}
Depending on the structures of $P_{2}$, explicit or numerical computation of the $D_{KL}$ between $P_{1}$ and $P_{2}$ may not be possible. In this case, the following proposition can be applied.
\begin{proposition}[Evidence Lower Bound]\label{prop:2}
Let $P_{2}, P_{3}\in Prob(\mathbb{R}^{n}\times \mathbb{R}^{k})$ absolutely continuous, and $P_{1}\in  Prob(\mathbb{R}^{n})$ absolutely continuous, let
\begin{equation}
ELBO(a,P_{2},P_{3})=\int\limits_{\mathbb{R}^{k}}\log(p_{2}(a|s))p_{3}(s|a)ds -D_{KL}(P_{3}(ds|a)||P_{2}(ds))
\end{equation}
Evidence Lower Bound, 
and let 
\begin{equation}
MELBO(P_{1},P_{2},P_{3})=\int\limits_{\mathbb{R}^{n}} (-\log(p_{1}(a))+ ELBO(a,P_{2},P_{3}))p_{1}(a)da
\end{equation}
the Mean Evidence Lower Bound.\\
It holds that 
\begin{equation}
-D_{KL}(P_{1}(da),P_{2}(da))\ge MELBO(P_{1},P_{2},P_{3})da.
\end{equation}
\end{proposition}
$P_{3}$ is called a variational distribution.\\
There exists a symmetric version of the KL divergence, the Jensen-Shannon divergence, which has the advantage of being defined for two general $P_{1},P_{2}\in Prob(\mathbb{R}^{n})$, without the absolutely continuity constraint.

\begin{definition}[Jensen-Shannon (JS) divergence] Let $P_{1},P_{2}$ probability measures on $\mathbb{R}^{n}\subseteq\mathbb{R}^{n}$, we define the Jensen Shannon divergence as  
\begin{equation}
D_{JS}(P_{1}||P_{2})=D_{KL}\left(P_{1}|| \frac{1}{2}P_{1}+\frac{1}{2}P_{2}\right)+D_{KL}\left(P_{2}|| \frac{1}{2}P_{1}+\frac{1}{2}P_{2}\right)
\end{equation}
\end{definition} This function is symmetric and is adopted in the training of the Classic Generative Adversarial Networks (the first ever introduced Generative Adversarial Network architecture) \cite{gan} and of Adversarial Autoencoders \cite{aae}, which are types of Generative Models.\\
In general, $D_{JS}$ cannot be computed analytically. However, lower bounds are available, as the following Proposition shows.\\
\begin{proposition} \label{prop:3}
Let $g:\mathbb{R}^{n}\rightarrow [0,1]$ continuous, $P_{1},P_{2}\in Prob({\mathbb{R}^{n}})$ then
\begin{equation}
2D_{JS}(P_{1}||P_{2})\ge 2\log2 + \int\limits_{\mathbb{R}^{n}}\log_{2}(g(a))P_{1}(da)+\int\limits_{\mathbb{R}^{n}}\log_{2}(1-g(a))P_{2}(da)   
\end{equation}
\end{proposition}
Using the formula above, and by applying Monte Carlo approximation, it is possible to approximate the JS divergence using samples from $P_{1}$ and $P_{2}$, without any other knowledge on the distributions.

Another important divergence is the Wasserstein distance:
\begin{definition}[Wasserstein Distance] 
Let $P_{1},P_{2}\in Prob(\mathbb{R}^{n})$, we define the Wasserstein distance as 
\begin{equation}
W_1(P_{1}, P_{2})=\inf_{\gamma \in \Gamma(P_{1}, P_{2})} \int_{\mathbb{R}^{n}\times \mathbb{R}^{n}} ||a-b||_{2}\gamma(da,db)\end{equation} where $\Gamma(P_{1}, P_{2})$  is the set of couplings of  $P_{1}$ and  $P_{2}$, i.e the set of probability distributions on $\mathbb{R}^{n}\times \mathbb{R}^{n}$ which have $P_{1}$ and $P_{2}$ as marginals. \end{definition} This function is symmetric and is adopted in Wasserstein Generative Adversarial Networks \cite{wgan}.

We remark that also the Wasserstein distance is generally hard to compute analytically. Because of the inverse triangular inequality, it holds
\begin{equation}\label{wlbound}
W_{1}(P_{1},P_{2})\ge |\int aP_{1}(da)-\int aP_{2}(da)|.
\end{equation}

We have the following proposition (which wraps theorems and propositions found in \cite{wgan,fineq,boundprob,ot}). \begin{proposition}
    
Let $P_{1}\in Prob(\mathbb{R}^{n})$ and $P_{m} \in Prob(\mathbb{R}^{n})$ $\forall n \in \mathbb{N}$.
Then
\begin{itemize}
    \item $\lim\limits_{m\rightarrow +\infty}\int\limits_{\mathbb{R}^{n}}f(a)P_{m}(da)=\int\limits_{\mathbb{R}^{n}}f(a)P_{1}(da) \forall f \in C(\mathbb{R}^{n})$ bounded or Lipschitz if and only if $\lim\limits_{m\rightarrow +\infty}W_{1}(P_{m},P)=0.$
    \item If $\lim\limits_{m\rightarrow +\infty}D_{JS}(P_{m},P)=0$ then $\lim\limits_{m\rightarrow +\infty}W_{1}(P_{m},P)=0$.
    \item If $\lim\limits_{m\rightarrow +\infty}D_{KL}(P_{1},P_{m})=0$ and $P_{m} \ll P_{1}$ $\forall m \in \mathbb{N}$ then $\lim\limits_{m\rightarrow +\infty}W_{1}(P_{m},P)=0$.
\end{itemize}

\end{proposition}

These properties justify the usage of the divergences defined above for learning Generative Models.

We remark that weak convergence of $P_{m}$ to  $P_{1}$ does not generally imply that $D_{KL}(P_{m},P_{1})$ converges to 0.  For example \cite{wgan} for the random variable $X_{m}=(\frac{1}{m},Z)$ with $Z \sim Uniform(0,1)$, which converges to $X=(0,Z)$, it holds that $D_{KL}$ is undefined $\forall m$, as $P_{1}$ and $P_{m}$ have disjoint supports, and $W_{1}(P_{1},P_{m})=\frac{1}{m}.$\\  Having established a foundational understanding of divergences, we now proceed to delineate Generative Models.
\subsection{Neural Networks and Generative Models}

Neural networks are, broadly speaking, statistical parametric models with a very high number of parameters. 
A type of neural networks are the Rectified Linear Unit (ReLU) neural networks \cite{gentheo}. \\
\begin{definition}[ReLU Feedforward Neural Network]
Let $L, N_0, \ldots, N_{L+1} \in \mathbb{N}$. A ReLU neural network with $L$ hidden layers is a collection of mapping $\phi: \mathbb{R}^{N_0} \mapsto \mathbb{R}^{N_{L+1}}$ of the form

\begin{equation} 
\phi(x)=A_L \circ \sigma \circ A_{L-1} \circ \cdots \circ \sigma \circ A_0(x) \quad x \in \mathbb{R}^{N_0},
\end{equation}
where $\sigma(x)=ReLU(x)$ is applied element-wise, $A_{\ell}(y)=M_{\ell} y+b_{\ell}$ is affine with $M_{\ell} \in \mathbb{R}^{N_{\ell+1} \times N_{\ell}}, b_{\ell} \in \mathbb{R}^{N_{\ell+1}}$, $\ell=0, \ldots, L$. The quantities $W=\max _{\ell=1, \ldots, L} N_{\ell}$ and $L$ are the width and depth of the neural network, respectively.
 We denote $\mathcal{N} \mathcal{N}(N_{0},N_{L+1},W, L)$ as the set of functions that can be represented by ReLU neural networks with width at most $W$ and depth at most $L$.\end{definition}

It can be proven \cite{gentheo} that a function $f_{*}(P_{2})$ (the push forward measure of $f$ through $\mu$), where $f:\mathbb{R}^{k}\rightarrow \mathbb{R}^{n}$ is a ReLU neural network and $P_{2}$ is an absolutely continuous measure on $\mathbb{R}^{k}$, can approximate with arbitrary precision any probability measure $P_{1}\in Prob(\mathbb{R}^{n})$, in the Wasserstein distance sense, even if $k<n$. \\

As being near in Wasserstein distance does not imply near similar in KL-divergence or in JS-divergence, we consider a Generative Model as 

\begin{equation}\label{eq:genmod}
p_{\alpha,\beta,\gamma}(a)=\int_{\mathbb{R}^{k}}p_{\gamma}(a-f_\alpha(s))p_{\beta}(s)ds
\end{equation}

where $f_{\alpha}:\mathbb{R}^{k} \rightarrow \mathbb{R}^{n}$ is a Neural Network of parameter $\alpha$, $P_{\beta}$ is an absolutely distribution on $\mathbb{R}^{k}$ with density $p_{\beta}$ and $P_{\gamma}$ is a distribution on $\mathbb{R}^{n}$. 
$P_{\gamma}$ is needed to deal with problems related to disjoint support caused by $f_{\alpha}$ having limited capacity \cite{survey}.
During inference, new samples will always be generated according to $(f_{\alpha})_{*}(P_{\beta})$.

In the following section, we provide an overview of the classes of Generative Models that we will adopt in the applications section.

\subsection{Examples of Generative Models}\label{subsec:gensec}

For ensuring reproducibility, we describe the algorithms of the Generative Models that we employ in the Applications section. Their content is by no means exhaustive. For a more detailed comparison, we refer to \cite{survey}.\\
With $P_{A}$, we denote the data generating distribution.
We introduce the concept of Autoencoder as it will be recalled many times in the following sections.
An Autoencoder is a function $f \circ g$ where $f:\mathbb{R}^{k}\rightarrow \mathbb{R}^{n}$, $g:\mathbb{R}^{n}\rightarrow \mathbb{R}^{k}$ are ReLU neural networks. In this framework, $f$ is called decoder while $g$ is called encoder. 

\subsubsection{Variational Autoencoder}

Variational Autoencoders have density

\begin{equation}
p_{\alpha,\gamma}(a,s)=p_{\gamma}(a-f_{\alpha}(s))p_{1}(s),
\end{equation}
where $P_{\gamma}(da)=\mathcal{MN}(0,\gamma I_{n})$, $P_{1}(ds)=\mathcal{MN}(0,I_{k})$ and $f_{\alpha}:\mathbb{R}^{k}\rightarrow \mathbb{R}^{n}$ is a decoder with weights $\alpha$, with $k<n$. $P_{\gamma},P_{1}$ are the distribution functions associated to the densities $p_{\gamma}$, $p_{1}$ respectively.
In order to train it, the log-likelihood function is approximated by the Mean Evidence Lower Bound (introduced in Proposition \ref{prop:2}) using a variational distribution $P_{\omega,\psi}(ds|a)=\mathcal{MN}(\mu_{\omega}(a),\sigma_{\psi}(a)^{2}I_{k})$ where $\mu_{\omega}:\mathbb{R}^{n}\rightarrow \mathbb{R}^{k}$ is an encoder with weights $\omega$,  and $\sigma_{\psi}:\mathbb{R}^{n}\rightarrow \mathbb{R}$ is a  neural network with weights $\psi$. With the symbol $\mathcal{MN}$ we denote a Multivariate Normal Distribution.\\ 
The associated optimization problem is
\begin{equation}
    \max_{\alpha,\omega,\psi,\gamma}MELBO(P_{A}(da),P_{\alpha, \gamma}(da,ds),P_{\omega,\psi}(da,ds)).
\end{equation}

The integrals are approximated using Monte Carlo integration.

\subsubsection{Adversarial Autoencoders}

Adversarial Autoencoders \cite{aae} have density $p_{\alpha}(a,s)=p_{2}(a-f_{\alpha}(s))p_{1}(s)$
where $P_{2}(da)=\mathcal{MN}(0, I_{n})$, $P_{1}(ds)=\mathcal{MN}(0,I_{k})$ and $f_{\alpha}:\mathbb{R}^{n}\rightarrow \mathbb{R}^{k}$ is decoder with weights $\alpha$, and $k<n$. $P_{2},P_{1}$ are the distribution functions associated with the densities $p_{2}$, $p_{1}$ respectively.

In this case a variational distribution $P_{\omega}(ds|a)=\delta_{\mu_{w}(a)}(ds)$ is adopted, where $\delta$ is the Dirac delta measure, $\mu_{\omega}:\mathbb{R}^{n}\rightarrow \mathbb{R}^{k}$ is an encoder with weights $\omega$, and the JS divergence is adopted instead of the KL divergence in the MELBO loss. The $D_{KL}$ is estimated according to Prop. \ref{prop:3} using a neural network $g_{\theta}:\mathbb{R}^{k}\rightarrow [0,1]$ with weights $\theta$. The optimization problem thus becomes \cite{aae}
\begin{eqnarray}
\max_{\theta}\min_{\alpha,\omega}\int\limits_{\mathbb{R}^{n}}\frac{1}{2}||f_{\alpha}(\mu_{w}(a))-a||_{2}^{2}P_{A}(da)+ \nonumber\\
+\int\limits_{\mathbb{R}^{n}}\left(\int\limits_{\mathbb{R}^{k}}\log_{2}(g_{\theta}(s))P_{\beta}(ds)+\int\limits_{\mathbb{R}^{k}}\log_{2}(1-g_{\theta}(s))P_{\omega}(ds|a)\right)P_{A}(da).
\end{eqnarray}
Also in this case, the integrals are approximated using Monte Carlo integration.

\subsubsection{Boundary Equilibrium Generative Adversarial Networks}
Boundary Equilibrium Generative Adversarial Networks \cite{began} assume that
$p_{\alpha}(a,s)=p_{2}(a-f_{\alpha}(s))p_{1}(s)$ where $P_{2}(da)=\mathcal{MN}(0,10^{-8}I_{n})$ and $p_{1}(s)=\mathcal{MN}(0,I_{k})$, $f_{\alpha}:\mathbb{R}^{k}\rightarrow \mathbb{R}^{n}$ is a decoder with weights $\alpha$ . In this case, no variational distribution is adopted and training is performed using an Autoencoder composed by an encoder $\mu_{\omega}:\mathbb{R}^{n}\rightarrow \mathbb{R}^{k}$ and a decoder $g_{\theta}:\mathbb{R}^{k}\rightarrow \mathbb{R}^{n}$ with weigths $\omega$ and $\theta$ respectively. Let $L_{\omega,\psi}(x)=||\mu_{\psi}(g_{\omega}(x))-x||_{2}$, the optimization problem is
\begin{equation}
    \min_{\theta,\omega}\max_{\alpha} \int\limits_{\mathbb{R}^{n}}L_{\omega,\theta}(a)P_{A}(da)-\int\limits_{\mathbb{R}^{k}}L_{\omega,\theta}(f_{\alpha}(s))P_{1}(ds),
\end{equation}
which is the lower bound of the Wasserstein distance (Eq. \ref{wlbound}) of the Autoencoder losses
while satisfying
\begin{equation}
    \int\limits_{\mathbb{R}^{n}}L_{\omega,\theta}(a)P_{A}(da)=\gamma \int\limits_{\mathbb{R}^{k}}L_{\omega,\theta}(f_{\alpha}(s))P_{1}(ds),
\end{equation}
with $\gamma>1$. This additional constraint is adopted to stabilise training, and in practice it is imposed using a closed-loop feedback control \cite{began}.
The integrals are approximated using Monte Carlo integration.

\subsubsection{Normalising Flows with Autoencoders}
Normalising flows models (NF) \cite{nf} when combined with an Autoencoder assume that $p_{\alpha,\beta}(a,s)=p_{2}(a-f_{\alpha}(s))p_{\beta}(s)$ where $P_{2}(da)=\mathcal{MN}(0, I_{n})$ and $p_{\beta}(s)$ is the density of $(g_{\beta})_{*}(P_{1})$, $f_{\alpha}:\mathbb{R}^{k}\rightarrow \mathbb{R}^{n}$ is a decoder with weights $\alpha$, $k<n$. The function $g_{\beta}:\mathbb{R}^{k}\rightarrow \mathbb{R}^{k}$ is $C^{2}$ neural network with a $C^{2}$ analytical inverse and weights $\beta$. Furthermore $P_{1}=\mathcal{MN}(0,I_{k})$. For training, a encoder $\mu_{\omega}:\mathbb{R}^{n}\rightarrow \mathbb{R}^{k}$, with weights $\omega$ is adopted. Due to the chain rule, the objective is
\begin{eqnarray}
    \min\limits_{\alpha,\omega,\beta}\int\limits_{\mathbb{R}^{n}} \left(||a-f_{\alpha}(\mu_{\omega}(a))||_{2}^{2}+\log p_{1}\left(g_{\beta}^{-1}(a)\right)\right)P_{A}(da) + \nonumber \\ 
     \int\limits_{\mathbb{R}^{n}} \log \left(\left|\det\left(\nabla\left[g_{\beta}^{-1}\right](a)\right)\right|\right)P_{A}(da).
\end{eqnarray}

\subsubsection{Energy Based Models with Autoencoders}

Energy Based Models (EBM) \cite{ebm} when combined with an Autoencoder assume that $p_{\alpha,\beta}(a,s)=p_{2}(a-f_{\alpha}(s))p_{\beta}(s)$ where $P_{2}(da)=\mathcal{MN}(0, I_{n})$ and $p_{\beta}(s)=\frac{e^{-g_{\beta}({s})}}{\int e^{-g_\beta({q})} d {q}}$, $f_{\alpha}:\mathbb{R}^{k}\rightarrow \mathbb{R}^{n}$, is a neural network with weights $\alpha$, $k<n$, and $g_{\beta}:\mathbb{R}^{k}\rightarrow \mathbb{R}$ is a neural network with weights $\beta$. For training, a encoder $\mu_{\omega}:\mathbb{R}^{n}\rightarrow \mathbb{R}^{k}$ with weights $\omega$ is adopted.

The objective is
\begin{equation}
\min_{\alpha,\omega,\beta}\int\limits_{\mathbb{R}^{n}} \left(||a-f_{\alpha}(\mu_{\omega}(a))||_{2}^{2}+\log p_{\beta}\left(\mu_{\omega}(a)\right)\right)P_{A}(da).
\end{equation}
However, in general, sampling from $p_{\beta}$ cannot be performed in an exact manner. It is known that $p_{\beta}$ is the stationary solution of the Ito Stochastic Differential Equation $ds_{t}=-\nabla g_{\beta}(s_{t})dt+dW_{t}$ where $W_{t}$ is a Brownian motion. As a consequence, sampling is possible by simulating different trajectories of the Stochastic Differential Equation until convergence.

\subsubsection{Denoising Diffusion Probabilistic models with Autoencoders}
For Denoising Diffusion Probabilistic Model (DDPM) \cite{ddpm}, the main idea is to slowly move from a known distribution to the data distribution in a large number of time steps while reversing a known process that has a fixed variance schedule $\xi_{t}\in [0,1], t=0,\dots, T$ using Variational Inference. When combined with Autoencoders, they assume that
$p_{\alpha,\beta}(a,s_{1})=p_{2}(a-f_{\alpha}(s_{1}))p_{\beta}(s_{1})$ where $P_{2}(da)=\mathcal{MN}(0, I_{n})$, $f_{\alpha}:\mathbb{R}^{k}\rightarrow \mathbb{R}^{n}$ is a decoder with weights $\alpha$, $k<n$, $P_{\beta}(ds_{1})$ is the marginal distribution of a distribution on $(\mathbb{R}^{k})^{T}$ such that
$P_{\beta}(s_{t-1}|s_{t})=MN\left(\frac{1}{\sqrt{\alpha_t}}\left({s}_t-\frac{\xi_t}{\sqrt{1-\bar{\eta}_t}} g_{\beta}({s}_t, t)\right), \xi_t I_n\right)$ $\forall t=2\ldots T$
where $g_{\beta}:\mathbb{R}^{k} \times \mathbb{R}\rightarrow \mathbb{R}$ is a neural network with weights $\beta$, $\eta_{t}=1-\xi_{t}$, $\bar{\alpha}_{t}=\prod_{s=1}^{t}\alpha_{s},$ $\kappa_{t}=\frac{\xi_t}{2 \eta_{t}\left(1-\bar{\eta}_t\right)}$ $
\forall t=1\ldots T$, and $P_{\beta}(s_{T})=MN(0,I_{k}).$
In this case a variational distribution $P_{2}\in Prob(\mathbb{R}^{k})^{T}$ such that
$P_{2}(s_{t}|s_{t-1})=MN(\sqrt{1-\xi_t} {s}_{t-1}, \xi_t I_n)$ is adopted. For training, a encoder $\mu_{\omega}:\mathbb{R}^{n}\rightarrow \mathbb{R}^{k}$ with weights $\omega$ is adopted, and the optimization problem is 

\begin{gather}
\min\limits_{\alpha,\omega,\beta}
\int\limits_{\mathbb{R}^{n}}(||a-f_{\alpha}(\mu_{\omega}(a))||_{2}^{2}P_{A}(da)+\\
\int\limits_{\mathbb{R}^{n}}\int\limits_{\mathbb{R}^{k}}\sum\limits_{t=1}^{T}\frac{1}{T}\kappa_{t}\left\|{s_{T}}-g_{\beta}\left(\sqrt{\bar{\alpha}_t} \mu_{\omega}(a)+\sqrt{1-\bar{\eta}_t} {s_{T}}, t\right)\right\|^2P_{\beta}(s_{T})P_{A}(da). \nonumber
\end{gather}\\

In the next section, we study how Generative Models can be employed to reduce the number of parameters of deformed geometries in order to improve DROMs and PPINNs for computational sciences.

\section{Application to Geometrically Parametrized Partial Differential Equations}\label{sec:two}
The remainder of this chapter will describe how Generative Models can be employed to improve the simulations of (possibly time-dependent) Geometrical Parametrized Partial Differential Equation (GPPDE) on $\mathbb{R}^{3}$ using DROMs or PPINNs approaches.
The most general form of a GPPDE is
\begin{equation}
\begin{cases}
\mathcal{L}(u,\gamma,t,x,y,z) & (x,y,z)\in \Omega(\gamma), \quad t\in (0,T], \quad \gamma \in \Gamma\\
\mathcal{B}(u,\gamma,t,x,y,z) & (x,y,z)\in \partial\Omega(\gamma), \quad t\in (0,T], \quad \gamma \in \Gamma\\
u(\gamma,0,x,y,z)=u_{0}(\gamma,x,y,z) & (x,y,z)\in \Omega(\gamma), \quad \gamma \in \Gamma
\end{cases}
\end{equation}
where $\gamma$ is a geometrical parameter, $\Gamma$ is the parameter space and $\mathcal{L}$ and $\mathcal{B}$ are operators acting on the interior $\Omega$ and boundary $\partial \Omega$ of the domain.\\
There are many methods to compute the value of the GPPDE solution for a certain value of $\gamma$ without explicitly solving the GPPDE.
One class of methods is the one of Non-Intrusive Reduced Order Model (DROM): discretise the parameter space $\Gamma$ and obtain $\gamma_{i}$, $i=1,\dots, m,$ then solve $m$ full order simulations on $N_{T}$ timesteps by discretising $\Omega(\gamma_{i})$ using classical techniques like Finite Difference, Finite Volume or Finite Element. If the discretised $\Omega(\gamma_{i})$ has the same number of points (which we assume to be equal $\frac{n}{3}$, where $n$ is a natural number multiple of 3) and the same topology. Once all simulations are performed, reduced order modelling could be applied to obtain a map \begin{equation}f: \Gamma \rightarrow \mathbb{R}^{N_{T}\times \frac{n}{3}}\end{equation} to reduce the computational cost of simulations with new parameters.

Another class of methods is the one of Parametrized Physics Informed Neural Network approach (PPINN): assume that $u$ is approximated by a neural network $u_{\theta}$, compute the loss function
    \begin{eqnarray}
    &\sum_{\gamma \in \Gamma}\sum_{(x,y,z)\in \Omega(\gamma)}\sum_{t\in (0,T]} \|\mathcal{L}(u_{\theta},t,x,y,z) \|^{2}+ \\
    &\sum_{\gamma \in \Gamma}\sum_{(x,y,z)\in \partial\Omega(\gamma)}\sum_{t\in (0,T]} \|\mathcal{B}(u_{\theta},t,x,y,z\|^{2}+ \\
    &\sum_{\gamma \in \Gamma}\sum_{(x,y,z)\in \Omega(\gamma)}||u_{\theta}(\gamma,0,x,y,z)-u_{0}(\gamma,x,y,z)||^{2}
    \end{eqnarray}
    and then optimize the loss. It is important to note that, because of how the loss is computed, each geometry need not follow the same discretization. Instead, a set of points in the domain, which can vary in size with $\gamma$ suffices.
    In this case, we directly learn a map
      \begin{equation}u_{\theta}: \Gamma\times [0, T]\times\mathbb{R}^3\rightarrow\mathbb{R}.\end{equation}\\

We will now see an application of Generative Models in solving GPPDEs and PPINNs.

\subsection{Application to Data Driven Reduced Order Modelling}

In this subsection, we study how to apply Generative Models when solving GPPDEs using classical numerical methods \cite{cgm}.
We assume to have $a_{1},..,a_{m}\in \mathbb{R}^{n}$ variables, each variable representing a discretised domain in $\mathbb{R}^{3}$. All underlying domains are assumed to have the same topology.
Note that if the $\gamma_{i}$ are not known the methodology is useful to find $s_{i}\in \mathbb{R}^{d}$ in order to parameterise all the point clouds. Instead, if they are known the methodology is useful to find the latent variables $s_{i}\in \mathbb{R}^{d}$ with $d<l$ in order to perform reduction in the parameter space. This chapter is focused on the latter.
The reduced or newly determined parameters can serve as inputs for solving GPPDEs.\\

Generative models can be integrated into a Data Driven Reduced Order Modelling framework by using the following pipeline:
\begin{itemize}
\item Generate constrained geometries ${a}_{1},\dots,{a}_{m}$ using Generative Models and store the corresponding latent variables ${s}_{1},\dots,{s}_{m}$, which will be used as a reduced parametrisation.
\item Preprocess the geometries (for example, using radial basis functions), e.g. to refine them.
\item Compute the full order solutions $u_{i}$ using Finite Difference/Finite Volume/Spectral Element/Finite Element methods.
\item Compute a non-intrusive reduced order map from the $s_{i}$ to the $u_{i}$.
\end{itemize}

\begin{svgraybox}
Generative models can also be employed to enforce a certain type of constraint on point clouds: multilinear constraints \cite{cgm}.
To define a multilinear constraint, we need to extract the $x$-components, $y$-components, and $z$-components of each point cloud. We consider a reference point cloud $a\in \mathbb{R}^{n}$.
By convention, we will assume that
\begin{itemize}
\item the $x$ coordinates of $a$ are the $\{a \cdot e_{j} | j\mod3 =0 \}$ where $e_{j}$ is the j-th basis element of $\mathbb{R}^{n},$ 
\item the $y$ coordinates of $a$ are the $\{a\cdot e_{j} | j\mod3 =1 \},$
\item the $z$ coordinates of $a$ are the $\{a\cdot e_{j} | j\mod3 =2 \}.$ 
\end{itemize}
We also assume that the coordinates are sorted such that $(x_{i},y_{i},z_{i})$ is the $i$-th point of the reference point cloud $\tilde{a}$ according to the mesh topology.
We recall that a  multilinear function is a function $(b_{1},\dots,b_{l}) :\mathbb{R}^{n}\rightarrow \mathbb{R}^{l}$ with variables $(x_{1},\dots,x_{\frac{n}{3}},y_{1},\dots,y_{\frac{n}{3}},z_{1},\dots,z_{\frac{n}{3}})$ such that:
\begin{itemize}
\item Each component $b_{j}$ is a polynomial with total order 3.
\item Each component $b_{j}$ has total order 1 when restricted to $(x_{1},\dots,x_{\frac{n}{3}})$, and therefore is affine when restricted to these variables.
\item Each component $b_{j}$ has total order 1 when restricted to $(y_{1},\dots,y_{\frac{n}{3}})$, and therefore is affine when restricted to these variables.
\item Each component $b_{j}$ has total order 1 when restricted to $(z_{1},\dots,z_{\frac{n}{3}})$, and therefore is affine when restricted to these variables.
\end{itemize}

Examples of this concept are as follows.
\begin{itemize}
\item $(x_{1}y_{2}z_{5}+x_{2}y_{3}+y_{5}+1,x_{1}+z_{1})$ is a multilinear function on two dimension,
\item $x_{1}x_{2}$ is not a multilinear function,
\end{itemize}
We define a  constraint of type $b(a_{i})=c\quad \forall i=1,\dots,m$ is multilinear.\\
To enforce this type of constraint, we adopt Algorithm \ref{alg:1}  \cite{cgm}.
An example of a multilinear function is the volume of a discretized surface triangular mesh with triangles with point index $\{ T_{h}=(i_{h},j_{h},k_{h} | h=1,\dots, N_{T}) \}$ where $N_{T}$ is the number of triangles
\begin{equation}
N_{T}=\frac{1}{6}\sum_{h=1}^{N_{T}} (-x_{k_{h}}y_{j_{h}}z_{i_{h}}+x_{j_{h}}y_{k_{h}}z_{i_{h}}+x_{k_{h}}y_{i_{h}}z_{j_{h}}-x_{i_{h}}y_{k_{h}}z_{j_{h}}-x_{j_{h}}y_{i_{h}}z_{k_{h}}+x_{i_{h}}y_{j_{h}}z_{k_{h}}).
\end{equation}
Note that steps 5), 8) and 9) can only be performed if the constraint is multilinear.

\end{svgraybox}

\begin{programcode}{Algorithm 2.1} \label{alg:1}
Input: $a\in \mathbb{R}^{n}$, a multilinear function $b$ and a constraint value $c$\\
Output: $\tilde{a}\in \mathbb{R}^{n}$ such that $b(\tilde{a})=c$ up to numerical errors.\\
Steps: \\
1) Compute $\delta c=c-b(a)$\\
2) Set $\tilde{a}=a$ and extract the point components $\tilde{x},\tilde{y},\tilde{z}$\\
3) Restrict $b$ to the $x$-coordinates and assemble the matrix $A_{yz}$ and the vector $b_{yz}$ such that $A_{yz}\tilde{x}+b_{yz}=b(\tilde{a})$\\
4) Solve the problem\\
\begin{equation}
\begin{aligned}
\min_{\delta x \in \mathbb{R}^{\frac{n}{3}}} \quad &  ||\delta x||^{2}\\
\textrm{s.t.} \quad & A_{yz}\delta x= \frac{\delta c}{3}\\
\end{aligned}
\end{equation}
which has solution
\begin{equation}
\delta x=(A_{yz}^{T}A_{yz})^{-1}A_{yz}\frac{\delta c}{3}
\end{equation}
5) Set $\tilde{x}=\tilde{x}+\delta x$\\
6) Restrict $b$ to the $y$-coordinates and assemble the matrix $A_{xz}$ and the vector $b_{xz}$ such that $A_{xz}\tilde{y}+b_{yz}=b(\tilde{a})$\\
7) Set 
\begin{equation}
\delta y=(A_{xz}^{T}A_{xz})^{-1}A_{xz}\frac{\delta c}{3}
\end{equation}
8) Set $\tilde{y}=\tilde{y}+\delta y$ \\
9) Restrict $b$ to the $z$-coordinates and assemble the matrix $A_{xy}$ and the vector $b_{xy}$ such that $A_{xy}\tilde{z}+b_{xy}=b(\tilde{a})$\\
10) Set 
\begin{equation}
\delta z=(A_{xy}^{T}A_{xy})^{-1}A_{xy}\frac{\delta c}{3}
\end{equation}
11) Set $\tilde{z}=\tilde{z}+\delta z$\\
12) Reassemble $\tilde{a}$ from $\tilde{x}$,$\tilde{y}$,$\tilde{z}$ \\
13) Return $\tilde{a}$.\\
\end{programcode}

The training geometries $a_{1},\dots,a_{m}$ were generated using Constrained Free Form Deformation \cite{cgm} (CFFD), implemented using the software PyGem \cite{pygem}.
To validate our pipeline, we compare the performances of different reduced order models (Gaussian Process Regression \cite{gpr}, Radial Basis Function Interpolation \cite{rbf}, Tree Regression \cite{tree}, implemented using the softwares EzyRB \cite{ezyrb} and Scikit-Learn \cite{sklearn}) on the original samples and those sampled from the Generative Model.\\
As a performance measure we adopt the following two errors:
\begin{equation}
max\_error=\max_{i=1\ldots m} \frac{|q_{i}-M_{-i}(\mu_{i})|}{|\max_{j=1\ldots m} q_{j}- \min_{j=1\ldots m}q_{j}|}
\end{equation}
and
\begin{equation}
mean\_error=\sum_{i=1\ldots m}\frac{1}{m} \frac{|q_{i}-M_{-i}(\mu_{i})|}{|\max_{j=1\ldots m} q_{j}- \min_{j=1\ldots m}q_{j}|},
\end{equation}

where $q_{i}$ is the $i$-th sample of the quantity of interest,$\mu_{i}$ is the corresponding geometric parameter, and $M_{-i}$ is the model trained on all data except $(\mu_{i},q_{i})$.
We chose this error measure because it is invariant by linear transformations of the output. In both test cases considered, $m=100$. We also compare the (averaged over 1000 runs) time of the computation of the error, which measures the average model speed.
We test the methodology using Variational Autoencoders (VAE), Adversarial Autoencoders (AAE), and Boundary Equilibrium Generative Adversarial Networks (BEGAN), described in the previous section. The output of $f_{\alpha}$ (defined as a neural network as in Section \ref{subsec:gensec}) is modified in order to impose the constraint of interest according to Algorithm \ref{alg:1}. \\
The algorithm is applied as the final layer of $f_{\alpha}$ both  during the training phase and the testing phase.

\begin{figure}
\centering
\includegraphics[width=\textwidth]{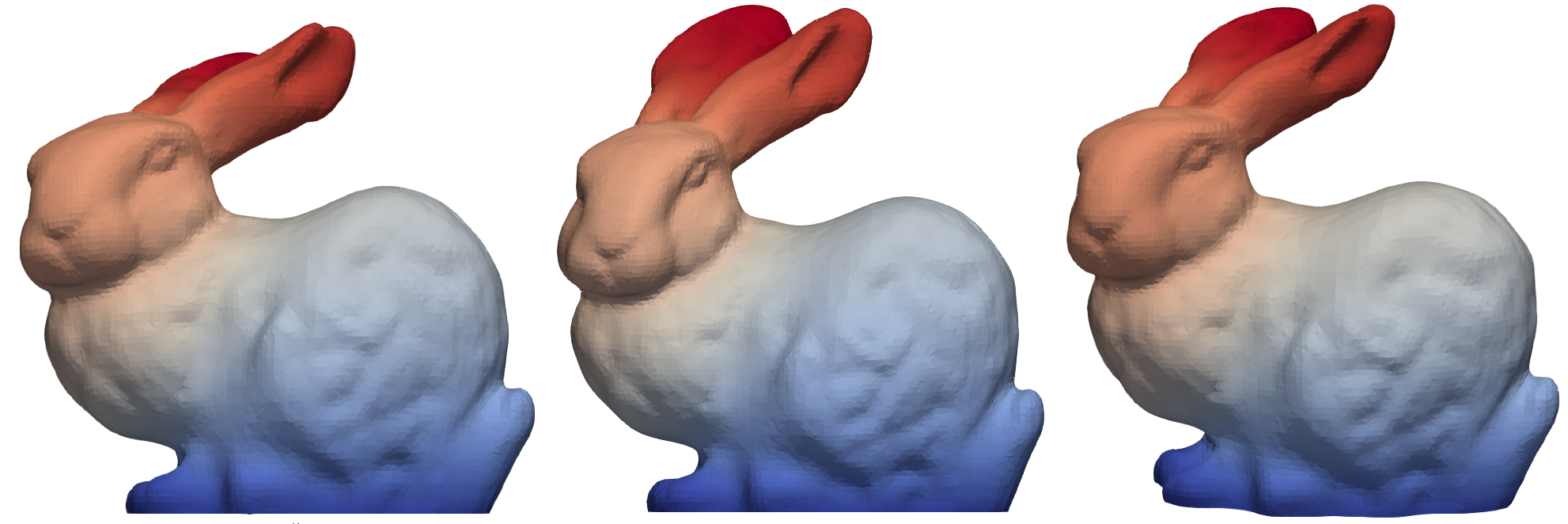}
\includegraphics[width=\textwidth]{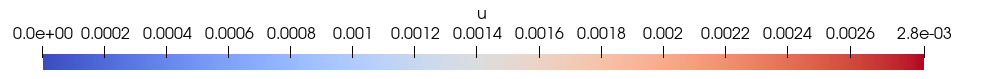}
\caption{Three simulations on meshes sampled using Generative Models.}
\label{fig:3}    
\end{figure}

\begin{figure}
\centering
\includegraphics[width=1\textwidth]{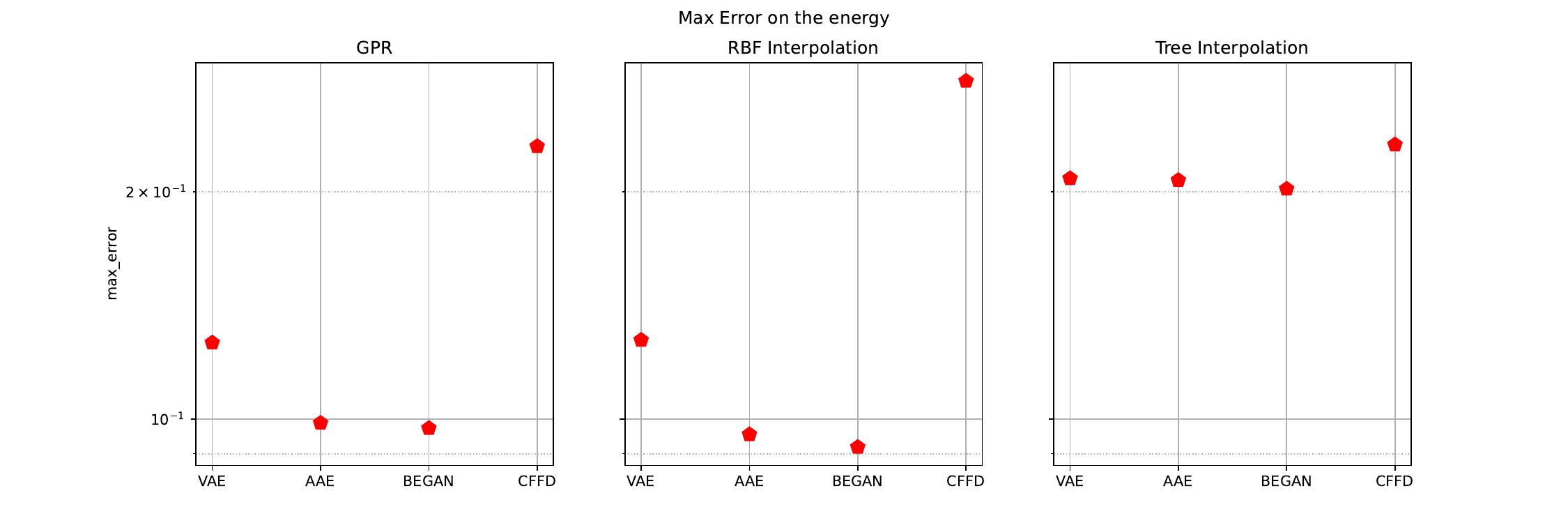}

\caption{Performance of different ROM techniques on 600 rabbits generated using CFFD and Generative Models in relation to the $max\_error$.}
\label{fig:4}  
\end{figure}

\begin{figure}
\centering
\includegraphics[width=1\textwidth]{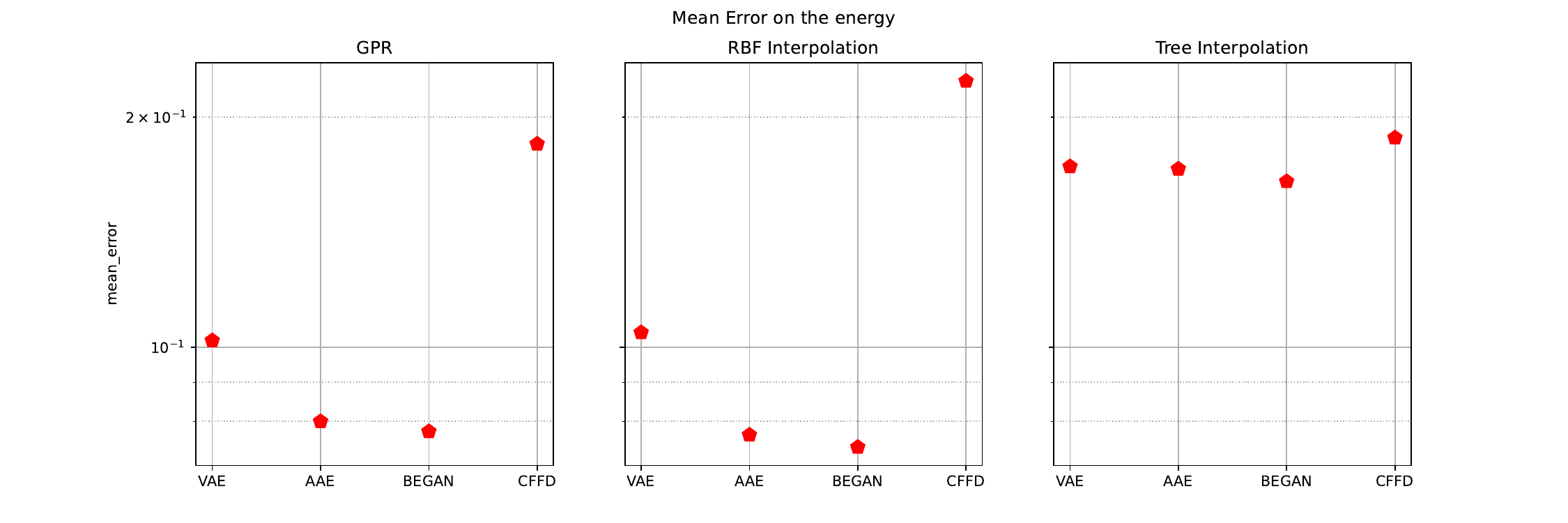}

\caption{Performance of different ROM techniques on 600 rabbits generated using CFFD and Generative Models in relation to the $mean\_error$.}
\label{fig:4_1}  
\end{figure}

\begin{figure}
    \includegraphics[width=1\textwidth]{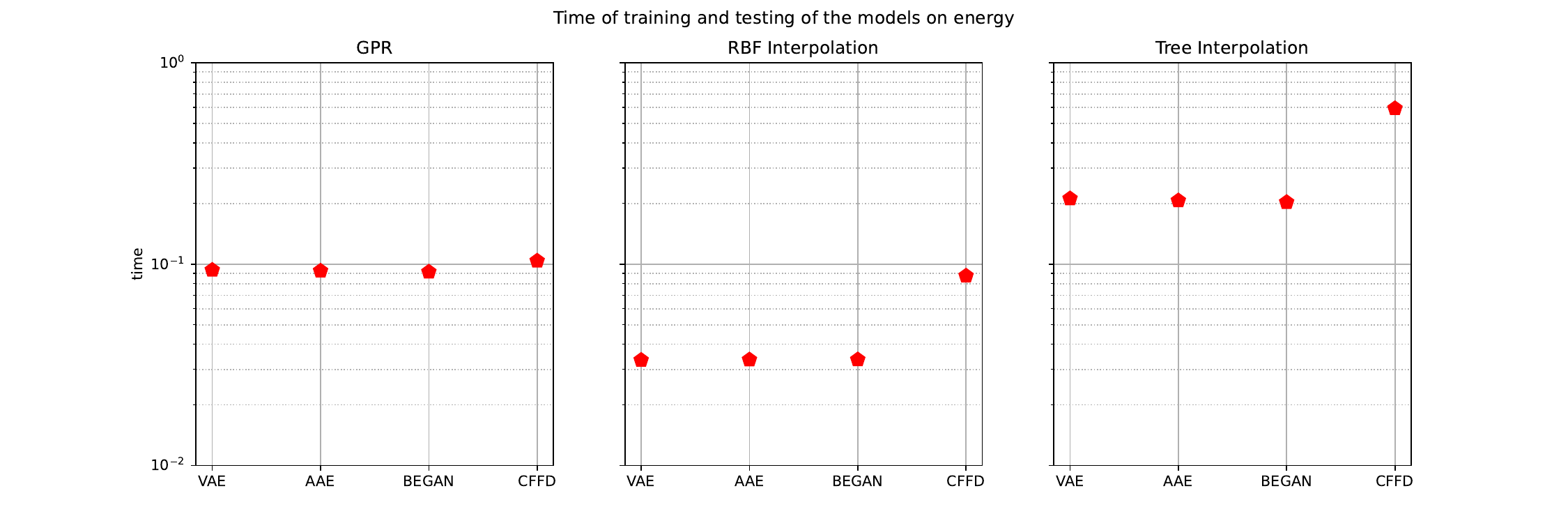}

\caption{CPU time for compute the $mean\_error$ averaged over 1000 simulations. The reduced parametrization also improves the time to compute the error.}
\label{fig:5}
\end{figure}

As a first example, we simulate a Poisson equation on a set of constant barycentre Stanford Bunny \cite{bunny} (see Fig. \ref{fig:3}) deformations :
\begin{equation}
\begin{cases}\nabla u(x,y,z)=f(x,y,z) & (x,y,z) \in \Omega \\ 
u(x,y,z)=0 & (x,y,z) \in \partial \Omega \cap\{y=0\} \\
\frac{\partial u}{\partial n}(x)=0 & (x,y,z) \in \partial \Omega \cap\{y\neq 0\}
\end{cases}
\end{equation}
where $\{y=0\}$ is the plane containing the base of the bunnies.
The dependency from $\gamma$ has been omitted for simplicity.\\
We are interested in the quantity $I_{u}=\int\limits_{\partial \Omega} u(x)dx$. 
 CFFD samples are obtained by moving control points following a uniform distribution in $[0,0.2]$, and then applying Algorithm \ref{alg:1}. A grid of control points of size $3\times 3 \times 3$ was adopted, and the control points near the base of the bunny were fixed, resulting in a total of 54 parameters.\\
Each bunny is made up of 145821 points in $\mathbb{R}^{3}$.\\
The barycenter formula is
\begin{equation}
\left(\frac{1}{n}\sum_{i=1}^{\frac{n}{3}}x_{i},\frac{1}{n}\sum_{i=1}^{\frac{n}{3}}y_{i},\frac{1}{n}\sum_{i=1}^{\frac{n}{3}}z_{i}\right).
\end{equation}
These samples compose the training data for the Generative Models.\\
As shown in Figs. \ref{fig:4}, \ref{fig:4_1}, \ref{fig:5} the reduced parametrisation obtained using Generative Models (each reduced parameter of a sample corresponds to the sample's latent representation in the Generative Model's latent space) increases the reduced order models performance in all the cases tested, Performance is increased as, due to the reduced parametrisation, the effects of  the curse of dimensionality are less severe (for example, Gaussian Process Regression generalisation error scales as $m^{-\frac{c}{dim({\Gamma})}}$, where $c$ is a constant independent of $m$ and $dim({\Gamma})$).  

As a second example, we simulate a two-phase water-air turbulent incompressible flow over a set of naval hulls with deformed bulbs (the original model is the Duisburg Test Case \cite{dtchull}, dependency from $s_{i}$ has been dropped for simplicity):
\begin{equation}
\begin{cases}
\partial_t(\rho {u})+\nabla \cdot(\rho {u} \otimes {u})+\nabla p-\rho g-\nabla \cdot \nu\nabla {u}-\nabla \cdot R=0    \\
\partial_t \alpha+\nabla \cdot({u} \alpha)=0\\
\alpha \rho_W+(1-\alpha) \rho_A =\rho \\
\alpha v\nu_W+(1-\alpha) \nu_A=\nu \\
\nabla \cdot {u}=0  
\end{cases}
\end{equation}

\begin{figure}
\includegraphics[width=0.5\textwidth]{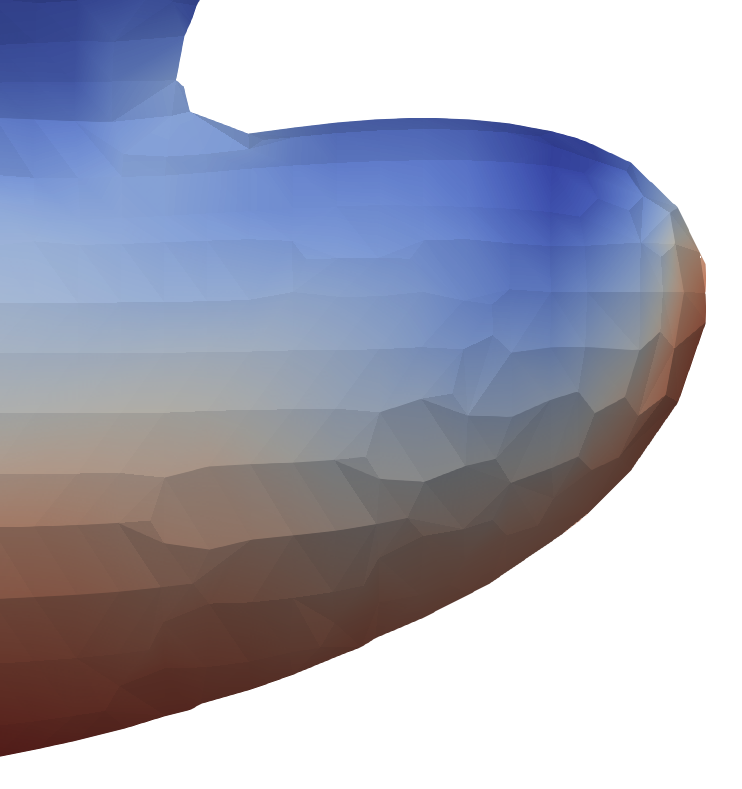}
\includegraphics[width=0.5\textwidth]{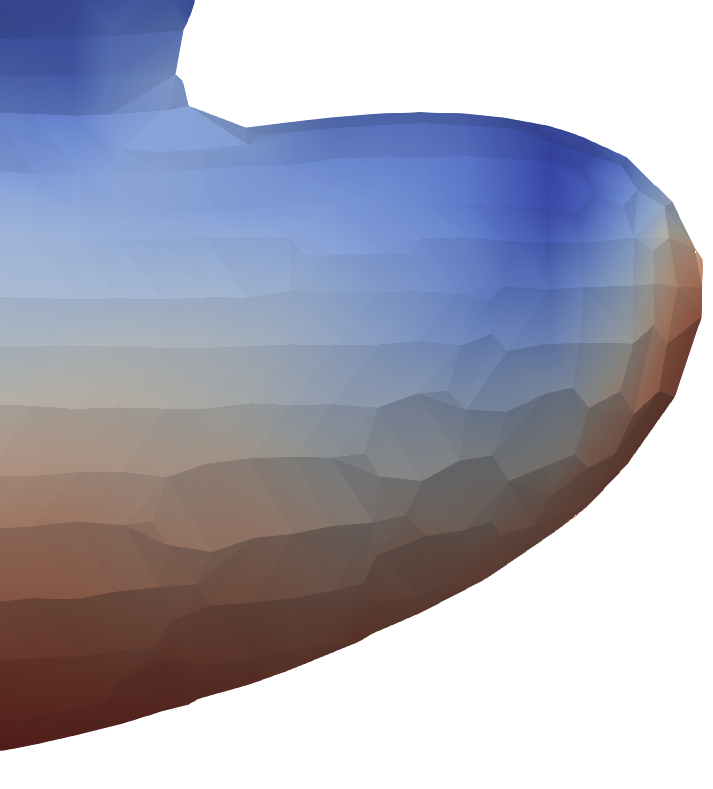}\\
\includegraphics[width=1\textwidth]{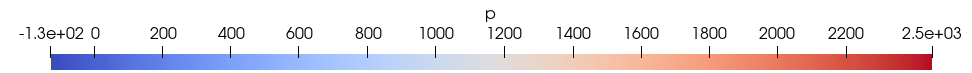}
\caption{Values of the pressure for simulations of two meshes sampled using Generative Models. Zoom only on the deformed part.}
\label{fig:6}     
\end{figure}

defined on 
\begin{equation}
(x,y,z,t) \in \Omega \times[0, T],
\end{equation}
where \begin{itemize}
\item $u$ is the velocity field, 
\item $p$ is the pressure field, 
\item $\rho_W, \rho_A$ are the densities of water and air, 
\item $v_W, v_A$ are the dynamic viscosities of water and air, 
\item $\rho$ is the mixed density,
\item $R$ is the Reynolds' stress tensor, 
\item $g$ is the acceleration of gravity,
\item$\alpha$ represents the interphase between the fluids with values from 0 (inside the air phase) to 1 (inside the water phase). 
\end{itemize}
Turbulence is modelled using the $\kappa-\omega$ model. For details on turbulence and boundary conditions, see \cite{cgm}.\\
We are interested in predicting the drag coefficient on the bulb \begin{equation}
c_{{d}}=\frac{1}{A_{h u l l}\left({u} \cdot {e}_x\right)^2}\left(\oint_{\delta \Omega_{\operatorname{man}}} p {n}-\left[\nu_\alpha\left(\nabla {u}+\nabla {u}^T\right)\right] {n} d s\right) \cdot {e}_x
\end{equation}
where $A_{hull}$ is the surface area of the hull.
Training geometries have been generated using CFFD, while full-order solutions are obtained using OpenFoam (see Fig. \ref{fig:6}).

\begin{figure}
\centering

\includegraphics[width=1\textwidth]{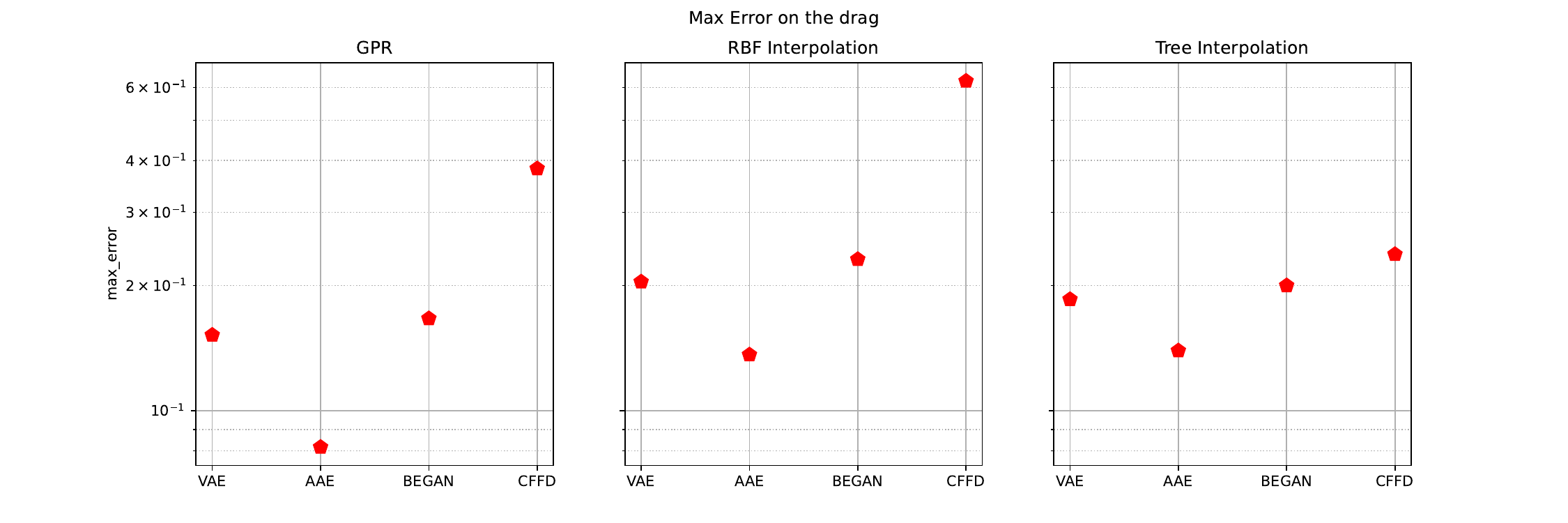}

\caption{Performance of different ROM techniques on 600 ships generated using CFFD and Generative Models in relation to the $max\_error$. }
\label{fig:7}
\end{figure}

\begin{figure}
\centering
\includegraphics[width=1\textwidth]{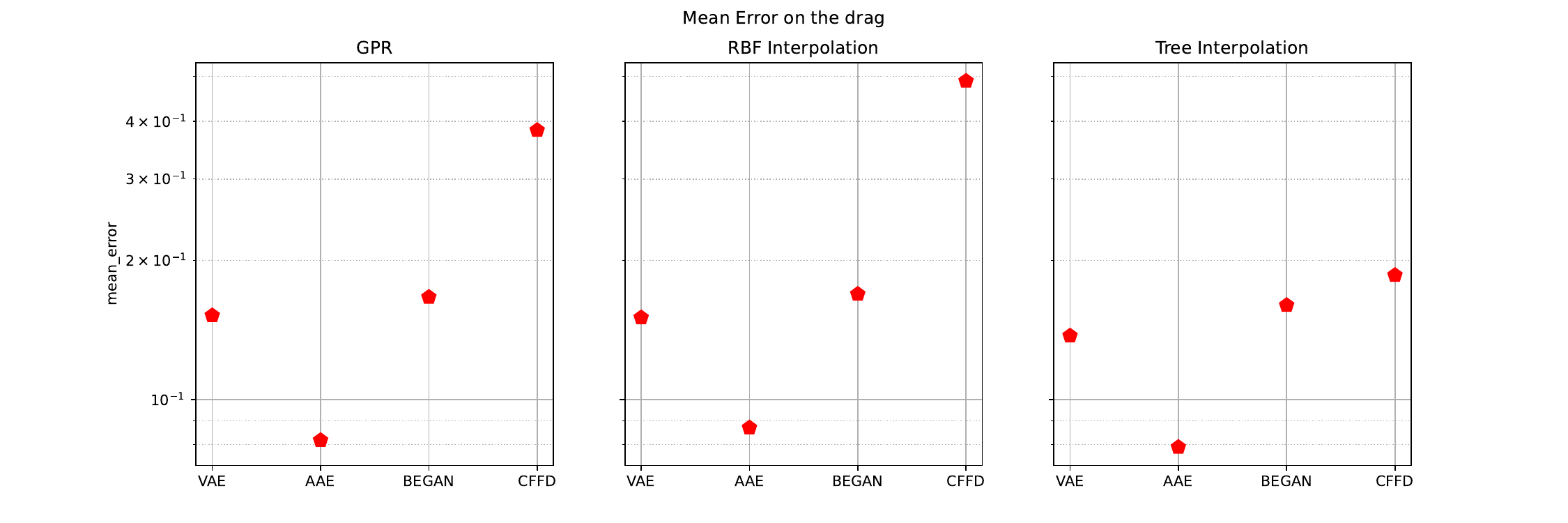}

\caption{Performance of different ROM techniques on 600 ships generated using CFFD and Generative Models in relation to the $mean\_error$. }
\label{fig:7_1} 
\end{figure}

\begin{figure}

\centering
\includegraphics[width=1\textwidth]{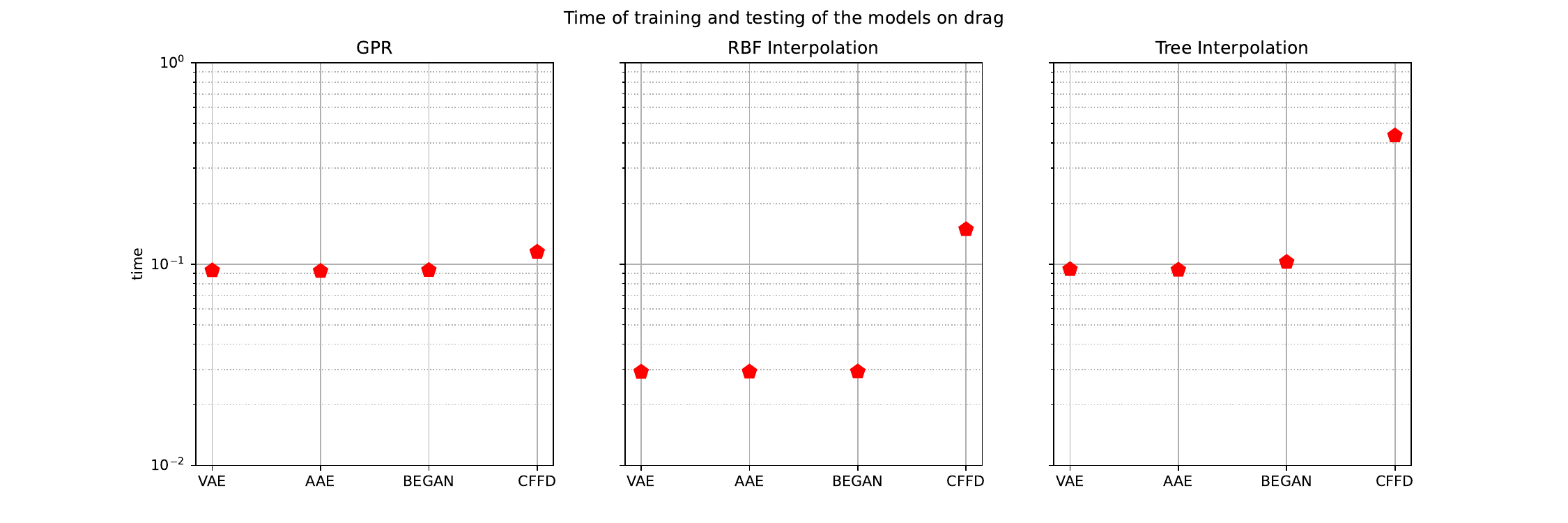}

\caption{CPU time for compute the $mean\_error$ averaged over 1000 simulations for the DTCHull test case.}
\label{fig:8}
\end{figure}

The CFFD geometries are obtained by moving the control points following a uniform distribution on $[0,0.2]$, and then by applying Algorithm \ref{alg:1}. A grid of control points of size $5\times 5 \times 5$ was adopted, and only the control points on the central subgrid $3\times 3 \times 3$ were allowed to change, resulting in $81$ parameters. \\
The reduced parametrisation (obtained using Generative Models) is 10.\\
Also in this more complex test case, Generative Models bring an increase in the performance of reduced order models (see Figs. \ref{fig:7}, \ref{fig:7_1}, \ref{fig:8}), as the effects of the curse of dimensionality are less severe.

\subsection{Point Conditioning and application to Parametrised Physics Informed Neural Networks}
The technique described in the previous subsection can be adapted to work with discretised domains with different topologies using a conditional autoencoder (inspired by PointNet \cite{pointnet}), which makes the technique applicable to PPINNs. 
The assumption is that all the deformed points cloud are in bijection with subsamples of the same point cloud $\bar{a}$. 
To keep the notation clean, we assume that all the points have the same topology, and we denote with ${point}_{j}(a_{i})$ the $j-th$ point of $a_{i}$.
The encoder is a map
\begin{equation}
g_{\psi}:\mathbb{R}^{3}\times \mathbb{R}^{3} \rightarrow \mathbb{R}^{d}
\end{equation} 
while the Decoder is a map
\begin{equation}
f_{\alpha}:\mathbb{R}^{d}\times \mathbb{R}^{3}  \rightarrow \mathbb{R}^{3}.
\end{equation}
The algorithm for computing the loss function is Algorithm \ref{alg:2}.
\begin{programcode}{Algorithm 2.2}
\label{alg:2}
Input: Samples $a_{1},\dots,a_{m}\in \mathbb{R}^{n}$ and a reference $\bar{a}\in \mathbb{R}^{n}$.
Output: Loss $l\in \mathbb{R}$.\\
1) $l=0$\\
2) for $i=1..m$:\\
3)          $\eta_{ij}=g_{\psi}(point_{j}(a_{i}),point_{j}(\bar{a}))$\\
4)          $s_{i}=\frac{3}{n}\sum_{j=1}^{n}\eta_{ij}$\\
5)	     $\tilde{a}_{i}=f_{\alpha}(point(\bar{a}),s_{i})$\\
6)          $l=l+||a_{i}-\tilde{a}_{i}||$\\
7) return $l$\\
\end{programcode}
The algorithm is employed to find a unique latent variable $s_{i}$ for each point cloud. The variable $s_{i}$ is then employed along a point $point_{j}(\bar{a})$ of the reference point cloud to obtain the corresponding point $point_{j}(\bar{a}_{i})$.
When the Autoencoder is trained, nonlatent variable models (such as Normalising Flows, Denoising Diffusion Probabilistic Models and Energy Based Models) can be trained to learn the distribution on the $s_{i}$, therefore we obtain a distribution $p_{A}(point_{j}(a)|point_{j}(\bar{a})) \text{  }\forall j$. This is a simplification of the approach proposed in \cite{pointflow}.
This algorithm, with respect to the methodology of the previous sections, supports $a_{i}$ with different topology and dimensions, as long as for each topology a discretization of the reference $\bar{a}$ of that dimension exists. Consequently $f_{\alpha}$ can be trained on a point cloud and then employed to deform a finer point cloud, which increases computational speed.
However, when sampling all points of a point cloud, the methodology is slower than the methods of the previous subsection, since $f_{\alpha}$ must be called $\frac{n}{3}$ times. \\
An application of this methodology is PPINN modelling, as in this case the loss is computed for one parameter $s_{i}$ and for one point $point_{j}(\bar{a}_{i})$ at a time.  As the total dimensionality of the input of the neural networks is reduced, fewer parameter and point samples are needed to perform correct inference. As a consequence, training is faster and the final loss is lower.\\
This methodology is similar to that of \cite{gapinn}, where the authors compute the latent variables of a set of non-parametric deformations employing a Variational Autoencoder and then leverage the latent variable to learn Parametric PINNs. However, authors of \cite{gapinn} adopt a Decoder like the one in the previous section: it is non conditioned and outputs only points clouds of the same topology. We relax this condition by allowing learning point clouds with different topologies. Furthermore, the authors do not explore the possibility of Generative Modelling, as the Variational Autoencoder is only employed to compute the latent representation.  We furthermore also study the efficacy of the approach with the perspective of parameter reduction.
We test the methodology on a Laplace equation on deformations of a Stanford Bunny obtained using Free Form Deformation \cite{ffd}. 
FFD samples are obtained by moving control points following a uniform distribution in $[0,0.2]$, and then applying the Algorithm \ref{alg:1}. A grid of control points of size $3\times 3 \times 3$ was adopted, and the control points near the base of the bunny were fixed, resulting in 54 parameters.\\
The latent parameterisation has dimension 5.
The equation is
\begin{equation}
\begin{cases}
\Delta u(x,y,z)=0 & (x,y,z)\in \Omega \\
u(x,y,z)=e^{z} & (x,y,z)\in \partial \Omega 
\end{cases}
\end{equation}

\begin{figure}

\includegraphics[width=0.5\textwidth]{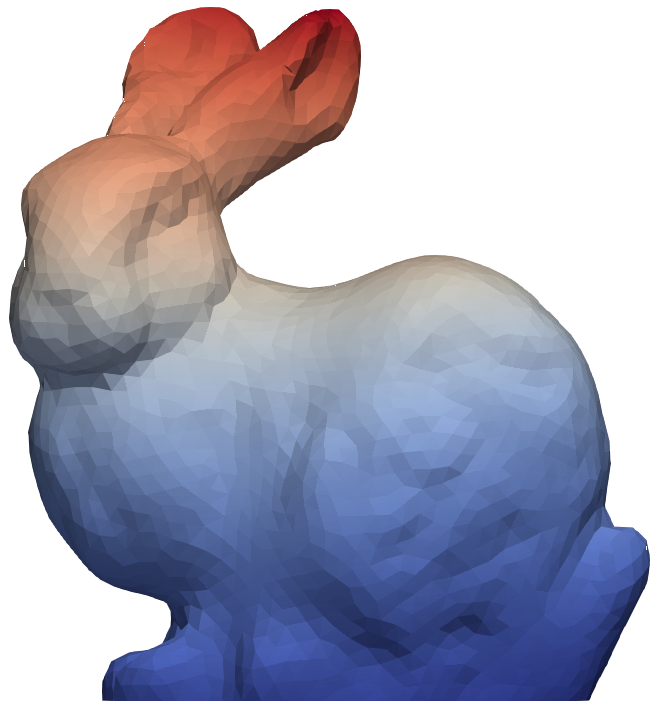}
\includegraphics[width=0.5\textwidth]{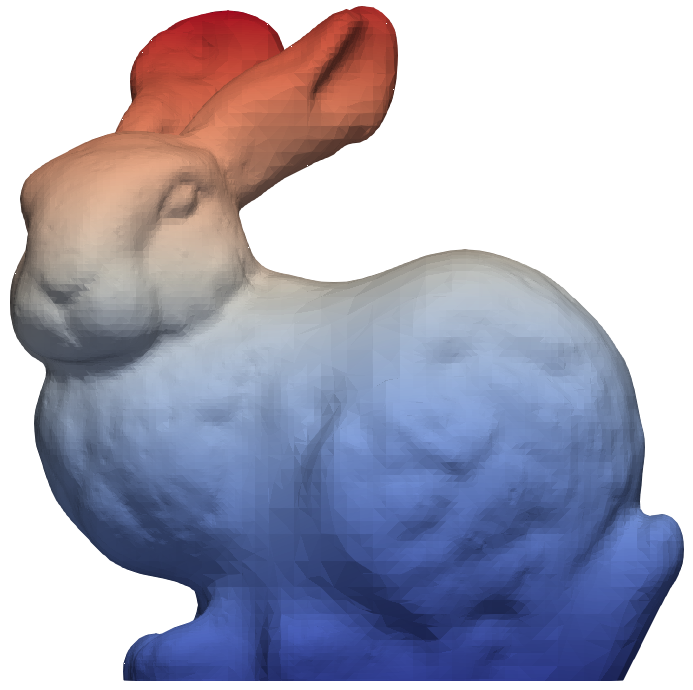}\\
\begin{center}
\includegraphics[width=1\textwidth]{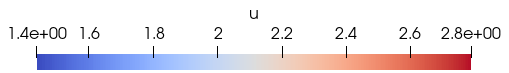}
\end{center}

\caption{The left figure is a bunny sampled with Generative Models that is employed as a training sample for the PPINN training. On the right a bunny sampled with higher resolution with respect to samples employed both for the Generative Model and for the PPINN training.}
\label{fig:9}   
\end{figure}

We test the methodology using Normalising Flows (NF), Denoising Diffusion Probabilistic Models (DDPM) and Energy Based Models (EBM).\\
The methodology can sample both the point clouds and the PDEs solution on point clouds with a higher resolution than the training ones (see Fig. \ref{fig:9}).
 Training has been done using PINA software \cite{pina}. 

\begin{figure}
\begin{center}
\includegraphics[width=0.45\textwidth]{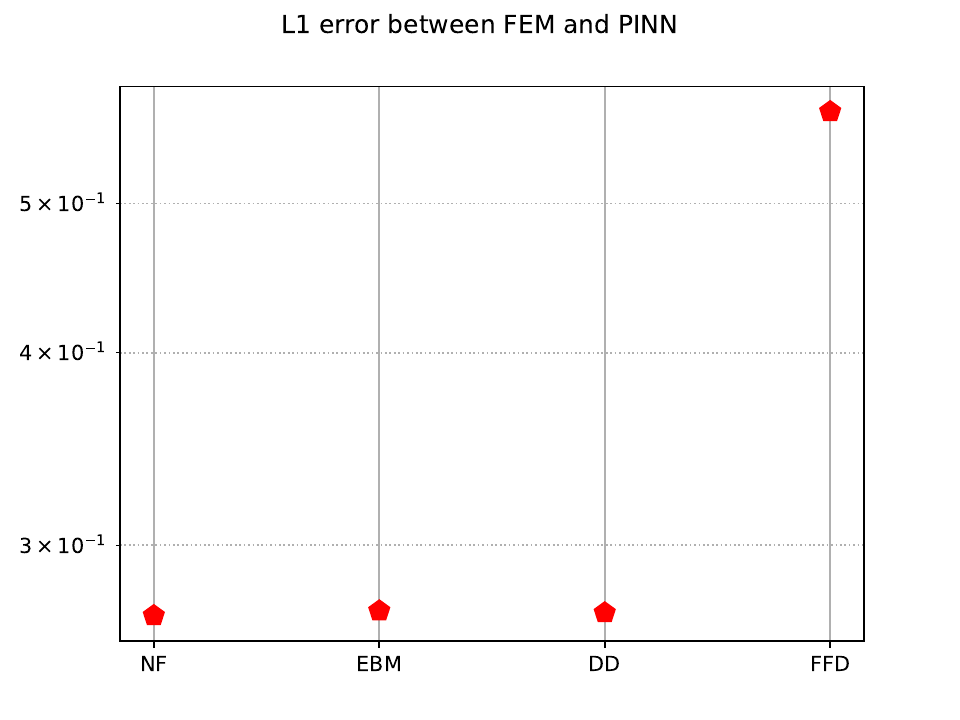}
\includegraphics[width=0.45\textwidth]{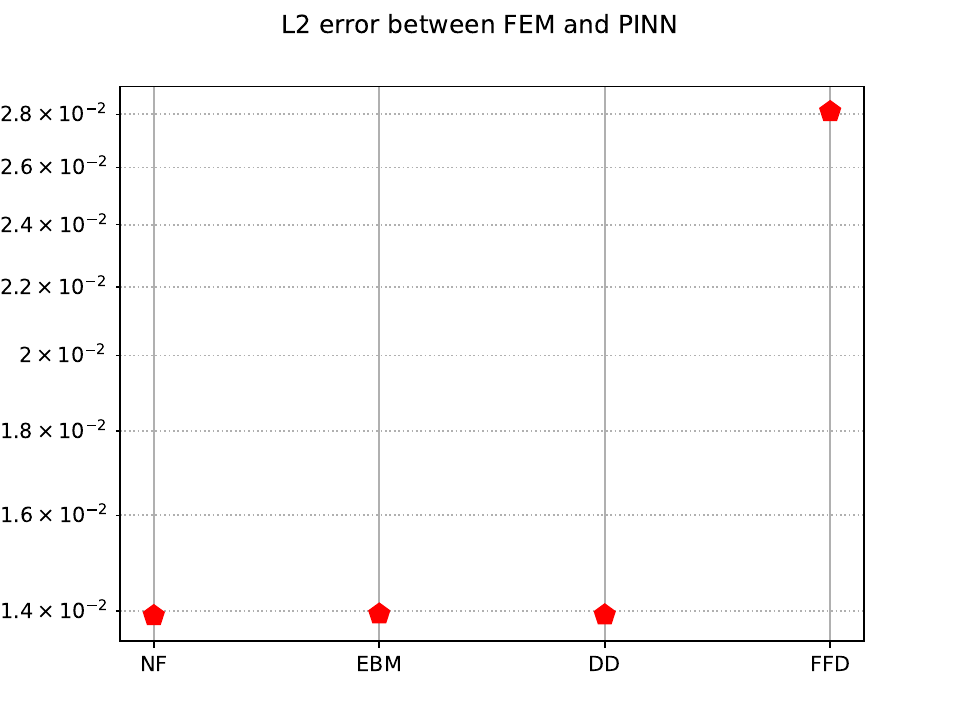}
\end{center}
\caption{L1 and L2 error between FEM solutions and respective PPINN ones on different parametrization}
\label{fig:10} 
\end{figure}

We validate our methodology by comparing the PPINN solution with the respective computed using the Finite Element Method (FEM). 
We adopt the following error metrics: 
\begin{equation}
L_{2}\_error=\frac{1}{m} \frac{1}{|\max_{j=1\ldots m} u_{j}- \min_{j=1\ldots m}u_{j}|}\sqrt{\sum_{i=1}^{m}||u_{PINN,i}-u_{FEM,i}||_{2}^{2}},
\end{equation}
and
\begin{equation}
L_{1}\_error=\frac{1}{m} \frac{1}{|\max_{j=1\ldots m} u_{j}- \min_{j=1\ldots m}u_{j}|}{\sum_{i=1}^{m}||u_{PINN,i}-u_{FEM,i}||_{1}},
\end{equation}
where $u_{FEM,i}$ represents the solutions obtained via the Finite Element Method on the coarse points, and $u_{PINN,i}$ are the respective PINN solutions. The results are shown in Fig. \ref{fig:10}. 

\begin{figure}
\begin{center}
\includegraphics[width=1\textwidth]{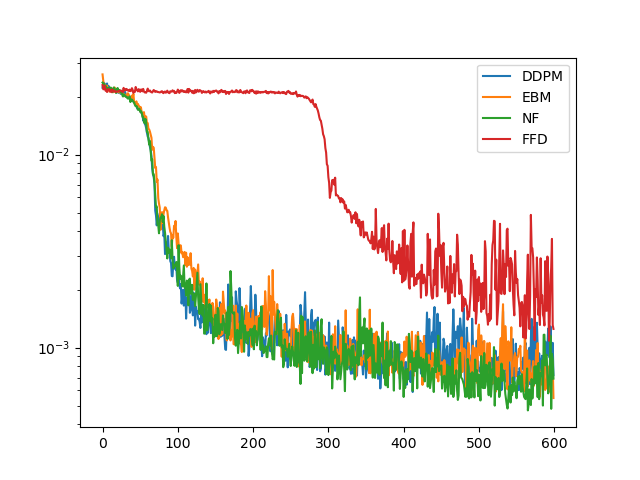}
\end{center}
\caption{Training loss of the PPINN of the laplace equation on geometries generated using different techniques. Overall, the training is better for geometries generated using Generative Models, as the number of parameters is lower. }
\label{fig:11}    
\end{figure}

Results are better when the parametrisation is obtained via Generative Models, as the reduced dimensionality simplifies the optimisation landscape, bringing a faster convergence and a lower local optima during the training process (Fig. \ref{fig:11}). In turn improves the L1 and L2 error.
Furthermore, the sampling is faster (by a 2x factor) thanks to the reduced input dimensionality.

\section{Conclusion}
\label{sec:3}
In this chapter, we introduced the foundational concepts of deep Generative Models, accompanied by practical examples. We demonstrated how Generative Models can effectively reduce the parameter space and enhance the performance of reduced order models, particularly for Geometrically Parametrised Partial Differential Equations (GPPDEs). We explored two methodologies: Data-Driven Reduced Order Models (DROMs), where we examined the feasibility of enforcing a multilinear constraint via Generative Models for parameter reduction, and we applied these to two test cases: a Poisson equation on a Stanford Bunny model and a real-world scenario involving multiphase flow on the DTCHull.
For Parameterized Physics-Informed Neural Networks (PPINNs), we developed a Generative Model-based technique to reduce the number of parameters, thereby accelerating convergence. This approach allows solutions to be learned on coarser meshes and then transferred to finer meshes, which we demonstrated through a Poisson equation.
Looking ahead, future work could involve testing Generative Models based on graph neural networks or applying Generative Models while preserving the data generation mechanism (Free From Deformation in this context) to ensure analytically sound sample quality.
\begin{acknowledgement}
We acknowledge the support provided by PRIN “FaReX - Full and Reduced order modelling of coupled systems: focus on non-matching methods and automatic learning” project, PNRR NGE iN-EST “Interconnected Nord-Est Innovation Ecosystem” project and INdAM-GNCS projects. This work was also partially supported by the U.S. National Science Foundation through Grant No. DMS-1953535 (PI A. Quaini).
\end{acknowledgement}

%
%
%

\appendix

\section{Appendix}
Here we describe the details of the architectures and hyperparameters adopted to train the Generative Models.\\
For what regards the BEGAN, AAE, VAE, each neural network has 6 hidden layers of size 500, with RELU activation functions, trained with  Batch Normalisation and Dropout. As an optimiser, AdamW the with learning rate of value $10^{-3}$ has been adopted.\\
For the Autoencoder adopted in combination with DDPM, EBM, and NF, we employ 4 hidden layers with depth 100, ReLU, and LayerNorm. As an optimiser, AdamW with 1e-03 has been adopted. Training is performed for 500 epochs.
For the neural network of the EBM and the DDPM, we employ a Feedforward with 6 hidden layers, of hidden dim 500 and ReLU. Training is done with Batch Normalisation and Dropout, is performed for 500 epochs. The optimized adopted is Adam with learning rate $10^{-3}$.  
For the NF, we employ 5 real-valued non-volume preserving layers. Training is performed for 500 epochs and Adam with learning rate $10^{-3}$.

\end{document}